\documentclass{article}
\usepackage{enumerate}
\usepackage{amsmath,amsfonts} 

\usepackage[applemac]{inputenc}

\newtheorem{theorem}{Theorem}[section]

\newtheorem{lemma}[theorem]{Lemma}
\newtheorem{definition}[theorem]{Definition}
\newtheorem{corollary}[theorem]{Corollary}

\newtheorem{remark}[theorem]{Remark}

\newcommand{\N}{{\mathbb N}}

\newcommand{\su}{ {\bf A}}
\newcommand{\Os}{\Omega_{\bf S}}

\begin{document}

\title{Measure-theoretic Uniformity and the Suslin Functional }
\author{Dag Normann\footnote{Department of Mathematics, The University 
of Oslo, P.O. Box 1053, Blindern N-0316 Oslo, Norway, email:  dnormann@math.uio.no}} 
\maketitle
\begin{abstract} \noindent We generalise results by Sacks and Tanaka concerning measure-theoretic uniformity for hyperarithmetical sets and a  basis theorem for $\Pi^1_1$-sets of positive measure to computability and semicomputability relative to the Suslin functional, alternatively to the (equivalent) Hyperjump.
\end{abstract}
\section{Introduction}
This note is a spinn-off of a project with Sam Sanders. In this project we combine methods from proof theory and higher computability theory to analyse the logical and computational complexity of classical theorems in analysis. These are theorems where the original formulations, referring to uncountable coverings and the like, are logically much more complex than the versions codeable in Second Order Arithmetic. 
The underlying problem inducing the research reported on in this note was about the relative computability of elements of classes of functionals (named $\Lambda$-functionals and $\Theta$-functionals, and introduced in \cite{paper1}) modulo the Suslin functional. The actual applications will appear in a paper still under planning. There the links between the Vitali covering theorem and $\Lambda$,  and between the Heine-Borel theorem and $\Theta$, will also be  discussed.

\medskip

The aim with this note is to lift the results by  Sacks \cite{Sacks.measure} and Tanaka \cite{Tanaka},  on measure-theoretic uniformity for hyperarithmetic theory, to computability relative to the Suslin functional.  The \emph{Suslin functional} is the functional of pure type 2 defined by
$${\bf S}(g) = \left \{ \begin{array}{ccc} 0 & {\rm if} & \forall f \in \N^\N \exists n \in \N (g(\bar f(n)) = 0) \\ 1 & {\rm if} & \exists f \in \N^\N \forall n \in \N (g(\bar f(n)) > 0)\end{array} \right.$$
The Suslin functional is computationally equivalent to the \emph{Hyperjump}, and is closely related to the $\Pi^1_1$-comprehension axiom. 

\medskip

It is well known that all subsets of $\N^\N$ computable, and even semicomputable, in ${\bf S}$ are $\Delta^1_2$, and since this is provable in ZFC, by \cite{FN} all such sets will be measurable for any completed Borel-measure on $\N^\N$ or on a closed subset thereof. So, the measurability of the sets we define will not be an issue.
\begin{remark}{\em  The actual meta-theorem says that if we can prove in ZFC that a set is $\Delta^1_2$ (with parameters), then we can also prove in ZFC that it is measurable. The essential forcing argument is due to R. Solovay \cite{Solovay}. Apparently, Solovay also observed (unpublished) the  meta-theorem  from \cite{FN}. } \end{remark}
 Some of the arguments in this note are well known from the classical literature on \emph{Descriptive Set Theory}, but the intention is to make the note reasonably self contained. The author's mental source is \cite{HJ}. For arguments on measure-theoretic uniformity, the reference will be Sacks \cite{Sacks}, and we refer to Kechris \cite{Ke} for some results in Descriptive Set Theory.

\medskip

In this note, we will first define the class of Suslin sets, an analogue of the Borel sets where we close under the Suslin operator instead of countable unions. (In the literature, the term \emph{Suslin set} is often used with a different definition, but we take the liberty to use our own concept here.) Then we show that all Suslin sets in the Cantor-space $C$ are measurable, and that we can compute the measure of a Suslin set from its code using ${\bf S}$. We then prove that the set of $f \in C$, relative to  which the least recursively inaccessible ordinal is the actual least recursively inaccessible, has measure 1. As a consequence we obtain that whenever $A \subseteq C$ is semicomputable in $\bf S$ and of positive measure, then $A$ contains an element computable in $\bf S$.
\subsection*{Acknowledgements} I am grateful to S. Sanders for involving me in the project that led to this research. I am also grateful to A.\ S.\ Kechris and G.\ E.\ Sacks for their brief comments on a preliminary note on the subject. Finally, I am grateful to the seminar on Mathematical Logic at the University of Oslo for comments during my presentation there.
\subsection*{For the preprint only}
The ease with which these results could be proved came as a surprise to the author. There are no new hard technical arguments in the note, they are just   adaptions of known methods to a situation needed by the author for an application. However, the results seem to be of general interest. 

If anyone reading this preprint, finds that the results overlap with those of others, please inform the author by e-mail.
\section{The Suslin Hierarchy}
\subsection{The Suslin sets}
Let SEQ be the set of finite sequences of integers. We use $s$, $t$ for elements of SEQ, and we identify them with the corresponding sequence numbers $\langle s \rangle$ and $\langle t \rangle$. We use the standard notation $\bar f(n) = \langle f(0) , \ldots , f(n-1)\rangle$.
\newline
We let $P$, $Q$ range over subsets of $X$, where $X$ may be the Cantor space $C$, the Baire space $\N^\N$ or any product of such spaces, possibly with  $\N$ as extra factors.

A \emph{Suslin scheme} $\bf P$ on $X$ is a map $s \mapsto P_s$, defined on SEQ. Given a Suslin scheme $\bf P$, we let $$\su({\bf P}) = \bigcup_{f \in \N^\N} \bigcap_{n \in \N} P_{\bar f(n)}.$$
$\su$ is known as the \emph{Suslin operator}.

By simultaneous recursion on the countable ordinal $\alpha$, we define two classes $\Sigma_\alpha$ and $\Pi_\alpha$ of subsets of $X$:
\begin{itemize}
\item[i)] $\Sigma_0$ consists of all sets that are both closed and open.
\item[ii)] For $\alpha > 0$, $\Sigma_\alpha$ is the set of $\su({\bf P})$ where $\bf P$ is a Suslin scheme with all $P_s$ in some $\Pi_\beta$, where $\beta < \alpha$ and may depend on $s$.
\item[iii)] $\Pi_{\alpha} = \{X \setminus P \mid P \in \Sigma_{\alpha}\}$.
\end{itemize}
\begin{definition}{\em Let $\Sigma_{\omega_1}$ be the union of all $\Sigma_\alpha$ for countable $\alpha $. We call the elements of $\Sigma_{\omega_1}$ for \emph{the Suslin sets}.}\end{definition}

It is easy to see that $\Sigma_\alpha \subseteq \Sigma_\gamma$ and $\Pi_\alpha \subseteq \Pi_\gamma$ when $\alpha \leq \gamma$, and that $\Pi_\alpha \subseteq \Sigma_{\alpha + 1}$ for all $\alpha$. Thus the complement of a set in $\Sigma_{\omega_1}$ will itself be in $\Sigma_{\omega_1}$. We also have
\begin{lemma}\label{lemma1} For each $\alpha > 0$, $\Sigma_{\alpha}$ is closed under the application of $\su$ to Suslin schemes with sets from $\Sigma_\alpha$, and consequently the class of Suslin sets is closed under the application of $\su$.\end{lemma} 
The proof is classical, and not too difficult. See e.g. Proposition (25.6)  in \cite{Ke}.

\medskip

In the same way we use well founded trees to code Borel sets, we can use well founded trees, or rather function codes for them, to code elements in $\Sigma_{\omega_1}$. We do not give any tedious details, it is easy to distinguish between codes for $\Sigma_0$, codes for $\Pi$-sets and codes for $\Sigma$-sets, such that when a tree codes a set, then we can extract codes for all the building blocks. The ordinal rank of the tree coded by a function $f$ will determine if $f$ is a $\sigma_\alpha$-code or a $\Pi_\alpha$-code.
The set of codes will be a complete $\Pi^1_1$-set.

\begin{remark}{\em If we consider the Borel sets as inductively defined from the clopen sets through the processes of countable unions and taking complements, each code for a Borel set can be viewed as a code for a Suslin set as well. Instead of interpreting a countable branching as representing a union, we interpret it via the use of the Suslin operator $\su$. However, we must distinguish between the two interpretations of a code. }\end{remark}

\begin{lemma} \label{lemma2}For each ordinal $\alpha$, the sets $\Sigma_\alpha$ and $\Pi_\alpha$ are closed under finite unions and intersections. Moreover, given a finite list of $\Sigma_\alpha$-codes or $\Pi_\alpha$-codes, we can compute the codes for the respective unions and intersections in the given classes uniformly in $^2E$.
\end{lemma}
{\em Proof}
\newline
By DeMorgan's laws it suffices to prove this for $\Sigma_\alpha$. Closure under finite unions is trivial. Closure under finite intersections requires a mild combinatoric effort, but essentially is proved  translating  the quantifier change
$$\forall k \leq m \exists f \forall n R(k,f,n) \mapsto \exists f \forall k \leq m \forall n R(k,(f)_k,n)$$ to a manipulation of Suslin schemes. This ends the proof.
\begin{definition}{\em A Suslin scheme ${\bf P}$ is \emph{normal} if $s \prec t \Rightarrow P_t \subseteq P_s$ for all $s$ and $t$ in SEQ.}\end{definition}
Every Suslin scheme is equivalent to a normal one in the sense that it defines the same set, see e.g. \cite{Ke},  and by Lemma \ref{lemma2} we can  reorganise a Suslin scheme to a normal one without increasing the rank. Using the recursion theorem for $^2E$, we may even transform a code for a well founded tree of Suslin schemes to a code for an equivalent tree of normal Suslin schemes, using a $^2E$-algorithm. From now on, we assume that all Suslin schemes are normal.

\medskip

For $\alpha > 0$ we can use the  arguments of Lemma \ref{lemma2} to  show that the classes of $\Sigma$-
sets and $\Pi$-sets are closed under countable unions and intersections, and that we  can find codes for the intersection and the union from a sequence of codes for such sets uniformly computable in $^2E$, but it is only for finite unions and intersections that we can preserve the rank.

\medskip

We end this section with a trivial application of the recursion theorem for $\bf S$:
\begin{lemma}
Each Suslin set is uniformly computable in $\bf S$ and any of its codes.\end{lemma}
Since the set of codes is computable in $\bf S$ as well, we see, by a standard diagonalization argument, that there will be sets computable in $\bf S$ that are not Suslin sets.
\subsection{Measurability} We will now restrict our attention to $X = C$, the \emph{Cantor space}.
Let $\bf m$ be the complete product measure on $C$. As stated in the introduction, all sets computable, or even semicomputable, in $\bf S$ will be measurable. In this section we will show that if $f$ is a code for the Suslin set $P[f]$, then we can compute ${\bf m}(P[f])$ from $f$ and $\bf S$ in the sense of computing its Dedekind cut.

\medskip

If $f,g \in \N^\N$, we let $f \leq g$ mean that $f$ is bounded by $g$ in the pointwise ordering. Likewise, if $s \in {\rm SEQ}$ and $f \in \N^\N$, we let $s \leq f$ mean that $s(i) \leq f(i)$ whenever $s(i)$ is defined.
\begin{definition}{\em \label{definition2.4}
Let $\bf P$ be a (normal) Suslin scheme and let $f \in \N^\N$. Let
$$\su_f({\bf P}) = \bigcup_{g \leq f} \bigcap_{n \in \N} P_{\bar g(n)}.$$

}\end{definition}
\begin{lemma}\label{lemma3}
With the notation from Definition \ref{definition2.4}, we have
$$\su_f({\bf P}) = \bigcap_{n \in \N} \bigcup_{lh(s) = n,s \leq f}P_s$$.
\end{lemma}
The proof is a trivial application of WKL. Here we need that the Suslin scheme is normal.

\medskip

In the next classical lemma we use that we are only dealing with sets that are measurable:
\begin{lemma}\label{lemma4} Let $\bf P$ be a Suslin scheme. Then
$${\bf m}(\su({\bf P})) = \sup \{{\bf m}(\su_f({\bf P})) \mid f \in \N^\N\}.$$
\end{lemma}
{\em Proof}
\newline
Let ${\bf m}(\su({\bf P})) = a$. Let $\epsilon > 0$. By recursion on $k$, we find $f(k)$ such that the union of all $\bigcap_{n \in \N}P_{\bar g(n)}$, where $g(i) \leq f(i)$ for $i \leq k$, has measure larger than $$a - \epsilon(1-2^{k+1}).$$ Then ${\bf m}(\su_f({\bf P})) \geq a - \epsilon$. This proves the lemma.
\begin{theorem}\label{theorem.2.10} Let $f$ be a code for the $\Sigma$-set or $\Pi$-set $P[f]$. Then the measure ${\bf m}(P[f])$ is computable in $\bf S$ uniformly in $f$.\end{theorem}
{\em Proof}
\newline
We use the recursion theorem, and argue by induction on the ordinal complexity of the code. The only non-trivial case is when $P[f] = \su({\bf P})$ where each $P_s$ is of the form $P[f_s]$. As an induction hypothesis we assume, with reference to Lemma \ref{lemma2}, that for any finite Boolean combination of the sets $P[f_s]$ we can compute the measure using $\bf S$. We can effectively enumerate all finite Boolean combinations of sets indexed by finite sequences $s$. Thus there will be a function $h = \langle h_a\rangle_{a \in \N}$, computable in $\bf S$ and $f$, such that $h_a$ gives us the measure of the corresponding Boolean combination of the sets $P[f_s]$. Using lemmas \ref{lemma3} and \ref{lemma4} we see that the measure of $\su(P[f])$, still seen as a Dedekind cut, will be $\Sigma^1_1$ in $h$, and thus computable in $\bf S$ and $f$ itself. This ends the proof.
\subsection{A basis theorem}
It is well known that any measurable set can be approximated from the inside and from the outside by Borel sets with the same measure. In this section we will show that we can find such sets for each $P[f]$ uniformly in $f$ and $\bf S$. We will use this to show that every Suslin set $A$ of positive measure contains an element computable in $\bf S$ and a code for $A$.

\medskip

\noindent The following should be well known:
\begin{lemma} \label{lemma5} Let $\bf P$ and $\bf Q$ be two Suslin schemes of measurable sets such that ${\bf m}(P_s \bigtriangleup Q_s) = 0$ for all finite sequences $s$, where $\bigtriangleup$ is the symmetric difference. Then ${\bf m}(\su({\bf P})) = {\bf m}(\su({\bf Q}))$. \end{lemma}
{\em Proof}
\newline
Since the measure of a Boolean  manipulation of countably many sets does not change through perturbations of the sets of measure 0, we can use Lemmas \ref{lemma3} and \ref{lemma4}.
\begin{theorem} \label{thm2.9} Let $f$ be a code for a $\Sigma$-set or $\Pi$-set $P[f]$. Then uniformly computable in $f$ and $\bf S$ we can find Borel codes for sets $B[f]$ and $C[f]$ such that $B[f] \subseteq P[f] \subseteq C[f]$ and ${\bf m}(B[f]) = {\bf m}(C [f])$. \end{theorem}
{\em Proof}
\newline
Once again, we use the recursion theorem and induction on the rank, and once again it is the case of an application of the Suslin operator $\su$ that is the only non-trivial one.
\newline
So, let $P[f]$ = $\su({\bf P})$ where $P_s = P[f_s]$ for each $s \in {\rm SEQ}$, and let $B[f_s]$ and $C[f_s]$ satisfy the induction hypothesis for each $s$. Then we can find a $g$ computable in $f$ and $\bf S$ such that $g$ is the join of the Borel-codes for all $B[f_s]$ and $C[f_s]$.
\newline
Now consider the Suslin schemes ${\bf B}$ and ${\bf C}$ given by the sets $B[f_s]$ and $C[f_s]$ respectively. Then both $\su({\bf B})$  and $\su({\bf C})$ are $\Sigma^1_1$-sets relative to $g$,  they have the same measure as $\su({\bf P})$ by Lemma \ref{lemma5}, and they are a subset, resp. superset of the same.
\newline
Let $a = {\bf m}(\su({\bf P}))$. For each $n \in \N$, we can use the Kleene basis theorem (see e.g. \cite{Sacks}, Theorem III 1.3) relative to $g$ to find $f_n$ such that ${\bf m}(\su_{f_n}({\bf B})) > a-2^n$. The sequence $\{f_n\}_{n \in \N}$ will be computable in $g$ and ${\bf S}$, and from this sequence and $g$ we can compute a Borel-code for 
$$\bigcup_{n \in \N} \su_{f_n}({\bf B}),$$  a set that will have measure $a$. We use this set as $B[f]$.
\newline
In order to find $C[f]$ we use the relativised version of the theorem saying that if a $\Pi^1_1$-set $D$ has a hyperarithmetical measure, there will be a $\Delta^1_1$-subset with the same measure (See \cite{Sacks}, Chapter IV). Thus, there will be a Borel-code hyperarithmetical in $g$ for a subset  $E$ of the set $C \setminus \su({\bf C})$ with the same measure as $C \setminus \su({\bf C})$ (Here $C$ is still the Cantor set). We use the complement of $E$ as our $C[f]$. This ends the proof.
\begin{corollary} Let $f$ be a code for a $\Sigma$-set or $\Pi$-set $P[f]$, and assume that ${\bf m}(P[f]) > 0$. Then $P[f]$ contains an element computable in $f$ and ${\bf S}$. \end{corollary}
This is an immediate consequence of Theorem \ref{thm2.9} and the Sacks-Tanaka basis theorem, see e.g. \cite{Sacks}.
\section{Computing with the Suslin-functional}

\noindent It is well known that if $F$ is a functional of type 2 computable in $^2E$, then the ordinal ranks of the computations will be bounded by a computable ordinal. ${\bf S}$ does not behave in a similar way. The set WO of codes for well orderings is computable in ${\bf S}$, and for each $f \in WO$, we can with ease, and in a uniform way, devise a terminating computation that lasts at least as long as the ordinal rank of $f$. 
However, we will see that the Suslin sets capture computations in ${\bf S}$ at a countable level much in the same way that the Borel sets capture computations relative to $^2E$.
\begin{definition}{\em Let $\Os$ be the set of sequences $\langle e , \vec a , \vec f , b\rangle$ such that $$\{e\}({\bf S},\vec a , \vec f) = b$$ and let $|| \cdot ||_{\Os}$ be the corresponding ordinal rank. (For notational reason, we assume that all arguments come in the order Suslin, integers, functions. $e$ is a \emph{Kleene-index} modified to this convention.)}\end{definition}
\begin{theorem}\label{theorem.main}
Let $\alpha$ be a countable ordinal.
\newline
For each $e$, $k$,  $\vec a$ and $b$ from $\N$,
$$\{\vec f \in (\N^\N)^k \mid ||\langle e, \vec a , \vec f , b\rangle ||_{\Os} < \alpha\}$$ will be a Suslin set.
\end{theorem}
{\em Proof}\newline
We prove this by induction on $\alpha$. In the case $\alpha = 0$, all sets in question are empty, so this case is trivial. When $\alpha$ is a limit ordinal, all sets in question are countable unions of Suslin sets by the induction hypothesis, so again, this case is trivial. We are then left with the \emph{simple} case when $\alpha = \beta + 1$, where we assume that the theorem holds for $\beta$. There will be 10 cases, the cases  S1 - S9 and \emph{otherwise}. For the cases S1 -  S3, S7, the initial computations, we do not even need the induction hypothesis, and for the \emph{otherwise}-case, all sets in question are empty. Case S9, the scheme of enumeration, and Case S6, the scheme of permutation,  follow by  direct applications of the induction hypothesis, since there is only one subcomputation in each case. Recursion, as defined in S5, is iterated composition, and the argument is as for S4. This leaves us with S4, composition, and S8, application of {\bf S}.
\begin{itemize}
\item[S4] Let $e = \langle 4 , e_1,e_2 \rangle$, and $$\{e\}({\bf S} , \vec a , \vec f) \simeq \{e_1\}({\bf S} , \{e_2\}({\bf S} , \vec a , \vec f),\vec a , \vec f).$$
Then we can use the induction hypothesis, and that the class of Suslin sets is closed under countable unions, since
$$\hspace*{-5mm}||\langle e , \vec a , \vec f , b\rangle||_{\Os} < \alpha \Leftrightarrow \exists c \in \N[||\langle e_2 , \vec a , \vec f , c\rangle||_{\Os} < \beta\; \wedge\; ||\langle e_1 , c,\vec a , \vec f , b\rangle||_{\Os}< \beta].$$
\item[S8] Application of ${\bf S}$: Let
$$\{e\}({\bf S}, \vec a , \vec f)  = \left\{ \begin{array}{ccc} 0& {\rm if}& \exists f \forall n \{e_1\}({\bf S} , \bar f(n) , \vec a , \vec f) = 0 \\1 & {\rm if} & \forall f \exists n \{e_1\}({\bf S} , \bar f(n) , \vec a , \vec f) > 0 \\
\bot&{\rm otherwise}& \end{array} \right.$$ where $\bot$ signifies that the value is undefined.

Let $R$ be the Suslin set
$$R = \bigcap_{s \in {\rm SEQ}} \bigcup_{c \in \N}\{\vec f \mid ||\langle e,s,\vec a , \vec f,c\rangle ||_{\Os} < \beta\}, $$ i.e. the set of $\vec f$ for which we apply $\bf S$ to a total function at level $\beta$. We will consider two Suslin schemes $\bf P$ and $\bf Q$ where
\begin{itemize}
\item $P_s  = R \cap \{\vec f \mid ||\langle e_1 , s , \vec a , \vec f,0 \rangle||_{\Os} < \beta\}$
\item $Q_s$ is  $R \setminus \bigcup_{c > 0}\{\vec f \mid ||\langle e_1 , s , \vec a , \vec f,c\rangle || < \beta\}$
\end{itemize}
Then $$ \su({\bf P})  = \{\vec f \mid ||\langle e , \vec a , \vec f, 0\rangle ||_{\Os} < \alpha\}$$ and $$R \setminus \su({\bf Q}) =  \{\vec f \mid ||\langle e,\vec a , \vec f , 1\rangle ||_{\Os} < \alpha\}.$$ In this case, for $b \neq 0,1$, the set in question will be empty, so we are trough.

\end{itemize}This ends the proof of the theorem.
\begin{remark}{\em If we have an ordinal code $h$ for $\alpha$ in the above proof, then a  code for the Suslin set constructed  will be computable in $h$ and $^2E$, uniformly in $h$. }\end{remark}
\subsection{Measure-theoretical uniformity for the Suslin functional}\label{subseq6.2}
We will now adjust arguments taken from \cite{Sacks}, and ``relativize" the result that almost all reals are hyperarithmetically low to $\bf S$. 

\medskip

The elements of the structure $L_{\omega_1^{\bf S}}[f]$ can be described as the interpretations $t^f$ of ranked terms $t = t({\bf f})$ where $\bf f$ is a constant that in each case is interpreted as $f$. The terms $t$ will be coded by elements in $L_{\omega_1^{\bf S}}$, and the interpretation $t^f$ will be in  $L_{\omega_1^{\bf S}}[f]$.
\newline
Any $\Delta_0$ formula for $L_{\omega_1^{\bf S}}$, with say  free variables x and X for numbers and sets,  will then be of the form
$$\phi(x,X,t_1, \ldots , t_k),$$
where $\phi$ has bounded quantifiers only. 
\newline
If $\phi(t_1, \ldots , t_k)$ is a $\Delta_0$-formula and $t_1 , \ldots , t_k \in L_\beta$ for $\beta < \omega_1^{\bf S}$, then $\{f \mid L_\beta[f] \models \phi\}$ will be a Suslin set with a code uniformly $^2E$-computable in codes for $\beta$ and $t_1 , \ldots , t_k$, and thus with a measure uniformly $\bf S$-computable in these codes.
\begin{lemma}
Let $0 \leq r \leq 1$, and let $\phi(x,X,t_1, \ldots , t_k)$ be a $\Delta_0$-formula where $x$ and $X$ are free variables. Assume that
$${\bf m}(\{f \mid \forall x \in \N \exists X \in L_{\omega_1^{\bf S}}[f] \phi(x,X,t_1^f , \ldots , t_k^f)\}) \geq r.$$
Then 
$${\bf m}(\{f \mid \exists \beta < \omega_1^{\bf S} \forall x \in \N \exists X \in L_\beta[f] \phi(x,X,t_1^f , \ldots , t_k^f)\}) \geq r.$$

\end{lemma}\label{6.3}
{\em Proof} 
\newline
It suffices to prove this for rational numbers $r$.
\newline 
We have that $$\{f \mid \forall x \in \N \exists X \in L_{\omega_1^{\bf S}}[f] \phi(x,X,t_1^f , \ldots , t_k^f)\}$$ $$ = \bigcap_{n \in \N}\{f \mid \forall x \leq n \exists X \in L_{\omega_1^{\bf S}}[f] \phi(x,X,t_1^f , \ldots , t_k^f)\}.$$
For each $n \in \N$ wee then have that 
$${\bf m}(\{f \mid \forall x \leq n \exists X \in L_{\omega_1^{\bf S}}[f] \phi(x,X,t_1^f , \ldots , t_k^f)\} )\geq r,$$
and for each rational $r' < r$ we can find $\beta_{n,r'} < \omega_1^{\bf S}$ such that 
$${\bf m}(\{f \mid \forall x \leq n \exists X \in L_\beta[f] \phi(x,X,t_1^f , \ldots , t_k^f)\} )> r'.$$ We  can use Gandy selection for $\bf S$ to find this $\beta_{n,r'}$. Since $r$ is rational, we can use the bounding principle in $L_{\omega_1^{\bf S}}$ and conclude that there is one $\beta < \omega_1^{\bf S}$ larger than all the $\beta_{n,r'}$. Then 
$${\bf m}(\{f \mid  \forall x \in \N \exists X \in L_\beta[f] \phi(x,X,t_1^t , \ldots , t_k^t)\}) \geq r.$$
This ends the proof of the theorem.
\begin{remark}{\em We have seemingly proved a stronger theorem than claimed, but the proof is only valid for the stronger claim when $r \in L_{\omega_1^{\bf S}}$.}\end{remark}

\begin{lemma}\label{6.4}
The set of $f \in C$ such that $\omega_1^{\bf S}$ is $f$-admissible has measure 1.
\end{lemma}
{\em Proof}
\newline
We only need to check $\Delta_0$-replacement, since the other Kripke-Platek axioms are trivially satisfied for all $f$. Since there are only countably many instances of replacement, we only need to verify the lemma for each one of them. By Lemma \ref{6.3}, the set of $f$ such that a particular instance of $\Delta_0$-replacement fails will have measure 0, so the set of $f$ for which this axiom holds will have measure 1.
\begin{lemma}\label{lemma.useful}\label{lemma6.1}
There is a well ordering $(A,\prec)$ of a subset of $\N$ of order type $\omega_1^{\bf S}$, semicomputable in $\bf S$, such that for each $a \in A$, $\{\langle b,c\rangle \mid b \prec c \prec  a\}$ is computable in $\bf S$, uniformly in $a$.

\end{lemma}
{\em Proof}
\newline
We let $A$ be the set of computation tuples $a = \langle e,\vec a , b\rangle$ such that $\{e\}({\bf S},\vec a) = b$ with the norm $|| \cdot ||$, and we let $a \prec a'$ if $||a|| < ||a'||$ or if $||a|| = ||a'||$ and $a < b$.
This ordering has the desired properties.
\begin{theorem}\label{6.5} The set of $f \in C$ such that $\omega_1^{\bf S} = \omega_1^{{\bf S},f}$  has measure 1. \end{theorem}
{\em Proof}
\newline
Let $(A,\prec)$ be as in Lemma \ref{lemma6.1},  let $a \in A$ and let $(A_a,\prec_a)$ be the corresponding initial segment of $(A,\prec)$, coded as the function $f_a$. 
\newline
Then applying a relativised version of a theorem due to Sacks \cite{Sacks.measure} and Tanaka \cite{Tanaka}, see also Corollary IV 1.6 of \cite{Sacks}, we have that $${\bf m}(\{f \mid \omega_1^{f_a,f} = \omega_1^{f_a}\}) = 1.$$

\noindent Taking the intersection of all these sets and the set from Lemma \ref{6.4}, we see that the set of $f$ such that $\omega_1^{\bf S}$ is both $f$-admissible and a limit of $f$-admissibles will have measure 1. This ends the proof.
\begin{corollary} Let $F:C \rightarrow \N$ be computable in ${\bf S}$. Then there is an ordinal $\alpha$ computable in ${\bf S}$ such that the algorithm for $F$ terminates at an ordinal rank below $\alpha$ on a set of measure 1. \end{corollary}
{\em Proof}
\newline
Let $F(f) = \{e\}({\bf S},f)$. By Theorem \ref{6.5}, the set of $f$ such that $\{e\}({\bf S},f)$  terminates before $\omega_1^{\bf S}$ has measure 1. Then for all $n \in \N$ there is an $\alpha_n < \omega_1^{\bf S}$ such that  the set of $f$ where $\{e\}({\bf S},f)$ terminates before $\alpha_n$ has measure $> 1-2^{-n}$. Combining the ${\bf S}$-effective version of Theorem \ref{theorem.main} and Theorem \ref{theorem.2.10} we can use Gandy selection for $\bf S$ to compute $\alpha_n$ from $n$. Thus $\{\alpha_n\}_{n \in \N}$ is bounded below $\omega_1^{\bf S}$, and we are through.
\begin{corollary} Let $A \subseteq C$ be semicomputable in ${\bf S}$ and with a positive measure.Then $A$ contains an element computable in $C$. \end{corollary}
{\em Proof}
\newline
Let $f \in A \Leftrightarrow \{e\}({\bf S},f)\!\!\downarrow$. By Theorem \ref{6.5}, almost all of these computations terminate at an ordinal level below $\omega_1^{\bf S}$. Then there is an $\alpha < \omega_1^{\bf S}$ such that the set of $f$ such that $\{e\}({\bf S},f)$ terminates before $\alpha$ will have positive measure.
\newline
This is a Suslin set of positive measure, and by the basis theorem for such sets, it contains an element computable in ${\bf S}$.

\end{document}